\documentclass[10pt, twoside]{article}
\usepackage{fancyhdr}
\usepackage{amssymb}
\usepackage{theorem}
\usepackage{array}
\usepackage{amsmath}
\usepackage{graphicx}

\newtheorem{theorem}{Theorem}[section]
\newtheorem{proposition}[theorem]{Proposition}
\newtheorem{lemma}[theorem]{Lemma}

\newtheorem{definition}[theorem]{Definition}
\newtheorem{remark}[theorem]{Remark}

\newtheorem{corollary}[theorem]{Corollary}
\newtheorem{conjecture}[theorem]{Conjecture}

\renewcommand{\arraystretch}{2}
\newlength{\blspace} 
\newsavebox{\qedbox} 

\sbox{\qedbox}{%
\begin{picture}(5.1,5)%
\put(0,0){\framebox(5,5){}}%
\end{picture}} %
\newcommand{\qed}{\hfill$\blacksquare$}%
\newenvironment{proof}[1][]%
{\begin{trivlist}\item{\textbf{Proof#1.}\hspace{\blspace}}}%
{\qed\end{trivlist}}

\fancypagestyle{plain}{%
   \fancyhead{} 
}

\newcommand{\gln}[1]{\mathrm{GL}_{#1}\mathbb{C}}
\newcommand{\bsln}[1]{\mathrm{SL}_{#1}\mathbb{C}}
\newcommand{\lgln}[1]{\mathfrak{gl}_{#1}\mathbb{C}}
\newcommand{\sln}[1]{\mathfrak{sl}_{#1}\mathbb{C}}

\renewcommand{\arraystretch}{1.0}

\begin{document}

\thispagestyle{plain}

\begin{center}
{\Large \textbf{A polynomiality property for Littlewood-Richardson coefficients}}
\bigskip

Etienne \textsc{Rassart}\\[2mm] 
\textit{Department of Mathematics, Massachusetts Institute of Technology}\\
\texttt{rassart@math.mit.edu}
\bigskip

\textrm{August 16, 2003}
\end{center}
\bigskip
       
\begin{abstract}
We present a polynomiality property of the Littlewood-Richardson coefficients $c_{\lambda\mu}^{\nu}$. The coefficients are shown to be given by polynomials in $\lambda$, $\mu$ and $\nu$ on the cones of the chamber complex of a vector partition function. We give bounds on the degree of the polynomials depending on the maximum allowed number of parts of the partitions $\lambda$, $\mu$ and $\nu$. We first express the Littlewood-Richardson coefficients as a vector partition function. We then define a hyperplane arrangement from Steinberg's formula, over whose regions the Littlewood-Richardson coefficients are given by polynomials, and relate this arrangement to the chamber complex of the partition function. As an easy consequence, we get a new proof of the fact that $c_{N\lambda\,N\mu}^{N\nu}$ is given by a polynomial in $N$, which partially establishes the conjecture of King, Tollu and Toumazet \cite{KTT} that $c_{N\lambda\,N\mu}^{N\nu}$ is a polynomial in $N$ with nonnegative rational coefficients.
\end{abstract}

\section{Introduction}

Littlewood-Richardson coefficients appear in many fields of mathematics. In combinatorics, they appear in the theory of symmetric functions (see \cite{Macdonald,EC2}). The Schur symmetric functions form a linear basis of the ring of symmetric functions, and the Littlewood-Richardson coefficients express the multiplication rule,
\begin{equation}
s_\lambda\cdot s_\mu = \sum_{\nu}c_{\lambda\mu}^{\nu}s_\nu\,,
\end{equation}
as well as how to write skew Schur function in terms of the Schur function basis:
\begin{equation}
s_{\nu/\lambda} = \sum_{\mu}c_{\lambda\mu}^{\nu}s_\mu\,.
\end{equation}

In the representation theory of the general and special linear groups, the characters of the irreducible polynomial representations of $\gln{k}$ are Schur functions in appropriate variables \cite{FultonHarris,Macdonald}. As such, the Littlewood-Richardson coefficient $c_{\lambda\mu}^{\nu}$ gives the multiplicity with which the irreducible representation $V_\nu$ of $\gln{k}$ appears in the tensor product of the irreducible representations $V_\lambda$ and $V_\mu$:
\begin{equation}
V_\lambda\otimes V_\mu = \bigoplus_{\nu}c_{\lambda\mu}^{\nu}V_\nu\,.
\end{equation}

Littlewood-Richardson coefficients also appear in algebraic geometry: Schubert classes form a linear basis of the cohomology ring of the Grassmannian, and the Littlewood-Richardson again express the multiplication rule \cite{Fulton}:
\begin{equation}
\sigma_\lambda\cdot \sigma_\mu = \sum_{\nu}c_{\lambda\mu}^{\nu}\sigma_\nu\,.
\end{equation}

In previous work with Billey and Guillemin \cite{BGR}, we studied the Kostka numbers $K_{\lambda\mu}$, which appear when expressing the Schur function $s_\lambda$ in terms of the monomial symmetric functions: $s_\lambda = \sum_{\mu}K_{\lambda\mu}m_\mu$. Kostka numbers also give the weight multiplicities in the weight space decomposition $V_\lambda = \bigoplus_{\mu}\big(V_\lambda\big)_{\mu}$ of the irreducible representation $V_\lambda$ of $\sln{k}$:
\begin{equation}
K_{\lambda\mu} = \mathrm{dim}\,\big(V_\lambda\big)_{\mu}\,.
\end{equation}  
We showed there that the Kostka numbers are given by a vector partition function and that this implies that the function $(\lambda,\mu)\longmapsto K_{\lambda\mu}$ is quasipolynomial in the cones of a chamber complex. We then defined a hyperplane arrangement, the Kostant arrangement, over whose regions this function was given by a polynomial. This allowed us to prove that the quasipolynomials in the cones were actually polynomials. As a corollary, we obtained an alternative proof to that of Kirillov that the function $N\mapsto K_{N\lambda\,N\mu}$ is a polynomial in $N$ for every fixed $\lambda$ and $\mu$.

In \cite{KTT}, King, Tollu and Toumazet conjecture that the Littlewood-Richardson coefficients exhibit a similar ``stretching'' property:

\begin{conjecture}{\rm\textbf{(King, Tollu, Toumazet \cite{KTT})}}
\label{conj:KTT}
For all partitions $\lambda$, $\mu$ and $\nu$ such that $c_{\lambda\mu}^{\nu}>0$ there exists a polynomial $P_{\lambda\mu}^{\nu}(N)$ in $N$ with nonnegative rational coefficients such that $P_{\lambda\mu}^{\nu}(0) = 1$ and $P_{\lambda\mu}^{\nu}(N) = c_{N\lambda\,N\mu}^{N\nu}$ for all positive integers $N$.
\end{conjecture}

In \cite{DerksenWeyman}, Derksen and Weyman prove the polynomiality part of this conjecture using semi-invariants of quivers. They call the functions $P_{\lambda\mu}^{\nu}(N)$ (for fixed $\lambda$, $\mu$ and $\nu$), \emph{Littlewood-Richardson polynomials}.

Here we extend the results of \cite{BGR} to the case of Littlewood-Richardson coefficients. We first express Littlewood-Richardson coefficients as a vector partition function (Theorem~\ref{thm:LRpartfunc}). This is done using a combinatorial model (the hive model \cite{Buch,KnutsonTao}) for computing the Littlewood-Richardson coefficients. This means that these coefficients are quasipolynomial in $\lambda$, $\mu$ and $\nu$ over the conical cells of a chamber complex $\mathcal{LR}_k$. 

From Steinberg's formula \cite{Steinberg}, giving the multiplicities with which irreducible representations appear in the decomposition into irreducibles of the tensor product of two irreducible representations of a complex semisimple Lie algebra, we then define a hyperplane arrangement, the Steinberg arrangement $\mathcal{SA}_k$. We show that the Littlewood-Richardson coefficients are given by a polynomial over the regions of this arrangement (Proposition~\ref{prop:SteinbergPolynomial}).

Finally, by comparing the chamber complex $\mathcal{LR}_k$ with the Steinberg arrangement $\mathcal{SA}_k$, we are able to show that the quasipolynomials in the cones of $\mathcal{LR}_k$ are actually polynomials in $\lambda$, $\mu$ and $\nu$, and we provide degree bounds (Theorem~\ref{thm:PolynomialityComplex}). Because we are working in cones, this provides an alternative proof to that of \cite{DerksenWeyman} of the polynomiality part of the conjecture of King, Tollu and Toumazet; we don't know whether the polynomials $P_{\lambda\mu}^{\nu}$ have nonegative coefficients or not. However, we get global polynomiality results in a chamber complex instead of polynomiality on fixed rays. We understand that Knutson \cite{Knutson} also proved polynomiality in cones using symplectic geometry techniques.

\subsection{Type $A$ root systems and Littlewood-Richardson coefficients}

The simple Lie algebra $\sln{k}$ (of type $A_{k-1}$) is the subalgebra of $\lgln{k} \cong \mathrm{End}(\mathbb{C}^k)$ consisting of traceless $k\times k$ matrices over $\mathbb{C}$. We will take as its Cartan subalgebra $\mathfrak{h}$ its subspace of traceless diagonal matrices. The roots and weights live in the dual $\mathfrak{h}^*$ of $\mathfrak{h}$, which can be identified with the subspace $x_1+\cdots+x_{k}=0$ of $\mathbb{R}^{k}$. The roots are $\{e_i-e_j\ : \ 1\leq i\neq j\leq k\}$, and we will choose the positive ones to be $\Delta_+ = \{e_i-e_j\ : \ 1\leq i<j\leq k\}$. The simple roots are then $\alpha_i = e_i-e_{i+1}$, for $1\leq i\leq k-1$, and for these simple roots, the fundamental weights are
\begin{equation}
\omega_i = \frac{1}{k}(\underbrace{k-i,k-i,\ldots,k-i}_{\textrm{$i$ times}},\underbrace{-i,-i,\ldots,-i}_{\textrm{$k-i$ times}})\,, \qquad 1\leq i\leq k-1\,.
\end{equation}
The fundamental weights are defined such that $\langle\alpha_i,\omega_j\rangle=\delta_{ij}$, where $\langle\cdot,\cdot\rangle$ is the usual dot product. The integral span of the simple roots and the fundamental weights are the root lattice $\Lambda_R$ and the weight lattice $\Lambda_W$ respectively. The root lattice is a finite index sublattice of the weight lattice, with index $k-1$. 

For our choice of positive roots, $\delta = \frac{1}{2}\sum_{\alpha\in\Delta_+}\alpha = \sum_{j=1}^{k-1}\omega_k = \frac{1}{2}(k-1,k-3,\ldots,-(k-3),-(k-1))$. The Weyl group for $\sln{k}$ is the symmetric group $\mathfrak{S}_k$ acting on $\{e_1,\ldots,e_k\}$ (i.e. $\sigma(e_i)=e_{\sigma(i)}$), and with the choice of positive roots we made, the fundamental Weyl chamber will be $C_0 = \{(\lambda_1,\ldots,\lambda_k)\ : \ \sum_{i=1}^k\lambda_i=0 \ \textrm{and } \lambda_1\geq\cdots\geq\lambda_k\}$. The action of the Weyl group preserves the root and weight lattices. Weights lying in the fundamental Weyl chamber are called \emph{dominant}, and we will call elements of the Weyl orbits of the fundamentals weights \emph{conjugates of fundamental weights}.

The finite dimensional representations of $\sln{k}$, or $\bsln{k}$, are indexed by the dominant weights $\Lambda_W\cap C_0$, and for a given dominant weight $\lambda$, there is a unique irreducible representation $\rho_\lambda : \sln{k} \rightarrow \mathfrak{gl}(V_\lambda)$ with highest weight $\lambda$, up to isomorphism. The finite dimensional polynomial representations of $\lgln{k}$, or $\gln{k}$, are indexed by partitions with at most $k$ parts, that is by sequences $(\lambda_1,\ldots,\lambda_k$) of integers satisfying $\lambda_1\geq\cdots\geq\lambda_k\geq 0$. Two irreducible representations $V_\lambda$ and $V_\mu$ of $\lgln{k}$ restrict to the same irreducible representation of $\sln{k}$ if $\lambda_i-\mu_i$ is some constant independent of $i$ for all $i$. So the irreducible representations of $\sln{k}$ correspond to equivalence classes of irreducible representations of $\lgln{k}$. Consider the map $\lambda\mapsto\bar{\lambda}$ given by
\begin{equation}
(\lambda_1,\ldots,\lambda_k) \quad\longmapsto\quad (\lambda_1,\ldots,\lambda_k) - \frac{\sum\lambda_i}{k}\underbrace{(1,1,\ldots,1)}_{\textrm{$k$ times}}\,.
\end{equation}  
Then the representations $V_\lambda$ of $\lgln{k}$ restricts to the irreducible representation $V_{\bar{\lambda}}$ of $\sln{k}$. Details about the construction of the irreducible representations of $\bsln{k}$ and $\gln{k}$ are well-known and can be found in \cite{Fulton} or \cite{FultonHarris}, for example. We will denote by $|\lambda|$ the sum $\sum\lambda_i$ (so $\lambda$ is a partition of the integer $|\lambda|$). We will also let $l(\lambda)$ denote the number of nonzero parts of $\lambda$.

Given two irreducible representations $V_\lambda$ and $V_\mu$ of $\gln{k}$, their tensor product $V_\lambda\otimes V_\mu$ is again a representation of $\gln{k}$, and we can decompose it in terms of irreducibles of $\gln{k}$:
\begin{equation}
\label{eqn:GLdecomp}
V_\lambda\otimes V_\mu = \bigoplus_{\nu}c_{\lambda\mu}^{\nu}V_\nu\,,
\end{equation}
where $c_{\lambda\mu}^{\nu}V_\nu = {V_\nu}^{\oplus c_{\lambda\mu}^{\nu}}$, for some nonnegative integer numbers $c_{\lambda\mu}^{\nu}$, called the \emph{Littlewood-Richardson coefficients}. The direct sum ranges over all partitions $\nu$, but $c_{\lambda\mu}^{\nu}=0$ unless $|\lambda|+|\mu|=|\nu|$ and $\lambda$ and $\mu$ are contained in $\nu$. We have a similar decomposition for the tensor product of two irreducible representations of $\sln{k}$:
\begin{equation}
\label{eqn:SLdecomp}
V_{\bar{\lambda}}\otimes V_{\bar{\mu}} = \bigoplus_{\bar{\nu}}m_{\bar{\lambda}\bar{\mu}}^{\bar{\nu}}V_{\bar{\nu}}\,,
\end{equation}
for nonnegative integers $m_{\bar{\lambda}\bar{\mu}}^{\bar{\nu}}$, where the sum ranges over all dominant weights $\bar{\nu}\in C_0$. 

There is a general formula due to Steinberg \cite{Humphreys,Steinberg} giving the multiplicity with which an irreducible representation $V_\nu$ occurs in the tensor product of two irreducible representations $V_\lambda$ and $V_\mu$ of a complex semisimple Lie algebra. This will give us a way of computing the $m_{\bar{\lambda}\bar{\mu}}^{\bar{\nu}}$, and also the $c_{\lambda\mu}^{\nu}$, but first we have to define the Kostant partitition function. 

\begin{definition}
The \emph{Kostant partition function} for a root system $\Delta$, given a choice of positive roots $\Delta_+$, is the function
\begin{equation}
K(v) = \Big|\Big\{(k_\alpha)_{\alpha\in\Delta_+}\in\mathbb{N}^{|\Delta_+|}\ : \ \sum_{\alpha\in\Delta_+}k_\alpha\alpha = v\Big\}\Big|\,,
\end{equation}
i.e. $K(v)$ is the number of ways that $v$ can be written as a sum of positive roots.
\end{definition}

\begin{theorem}{\rm \textbf{(Steinberg \cite{Steinberg})}}
\begin{equation}
\label{eqn:Steinberg}
m_{\bar{\lambda}\bar{\mu}}^{\bar{\nu}} = \sum_{\sigma\in\mathfrak{S}_k}\sum_{\tau\in\mathfrak{S}_k}(-1)^{\mathrm{inv}(\sigma\tau)}K(\sigma(\bar{\lambda}+\delta)+\tau(\bar{\mu}+\delta)-(\bar{\nu}+2\delta))\,,
\end{equation}
where $\mathrm{inv}(\psi)$ is the number of inversions of the permutation $\psi$.
\end{theorem}

Restricting equation~\eqref{eqn:GLdecomp} to $\bsln{k}$, we get
\begin{equation}
V_{\bar{\lambda}}\otimes V_{\bar{\mu}} = \sum_{\nu}c_{\lambda\mu}^{\nu}V_{\bar{\nu}}\,,
\end{equation}
and comparing with~\eqref{eqn:SLdecomp} gives
\begin{equation}
\label{eqn:SLandGL}
c_{\lambda\mu}^{\nu} = m_{\bar{\lambda}\bar{\mu}}^{\bar{\nu}}\,.
\end{equation}
Hence Steinberg's formula also computes the Littlewood-Richardson coefficients, and we can further simplify things by noticing that if we let $\mathbf{1}_k$ denote the vector $(1,1,\ldots,1)\in\mathbb{R}^k$, then
\begin{eqnarray*}
\sigma(\bar{\lambda}+\delta)+\tau(\bar{\mu}+\delta)-(\bar{\nu}+2\delta) & = & \sigma(\bar{\lambda})+\tau(\bar{\mu}) - \bar{\nu} + \sigma(\delta) + \tau(\delta) - 2\delta\\
& = & \sigma(\lambda-\frac{|\lambda|}{k}\mathbf{1}_k) + \tau(\mu-\frac{|\mu|}{k}\mathbf{1}_k) - (\nu-\frac{|\nu|}{k}\mathbf{1}_k) + \sigma(\delta) + \tau(\delta) - 2\delta\\
& = & \sigma(\lambda)-\frac{|\lambda|}{k}\mathbf{1}_k + \tau(\mu)-\frac{|\mu|}{k}\mathbf{1}_k - \nu+\frac{|\nu|}{k}\mathbf{1}_k + \sigma(\delta) + \tau(\delta) - 2\delta\\
& = & \sigma(\lambda+\delta) + \tau(\mu+\delta) - (\nu+2\delta) + \frac{1}{k}(|\nu|-|\lambda|-|\mu|)\mathbf{1}_k\\
& = & \sigma(\lambda+\delta) + \tau(\mu+\delta) - (\nu+2\delta)\,.
\end{eqnarray*}

In view of~\eqref{eqn:Steinberg} and~\eqref{eqn:SLandGL}, this gives
\begin{equation}
c_{\lambda\mu}^{\nu} = \sum_{\sigma\in\mathfrak{S}_k}\sum_{\tau\in\mathfrak{S}_k}(-1)^{\mathrm{inv}(\sigma\tau)}K(\sigma(\lambda+\delta)+\tau(\mu+\delta)-(\nu+2\delta))\,.
\end{equation}

In Section~\ref{sec:SteinbergArrangement}, we will use this formula to define a hyperplane arrangement over whose regions the Littlewood-Richardson coefficients are given by polynomials in $\lambda$, $\mu$ and $\nu$.

\subsection{Partition functions and chamber complexes}
\label{subsec:ChamberComplexes}

Partition functions arise in the representation theory of the semisimple Lie algebras in the form of Kostant's partition function, which sends a vector in the root lattice to the number of ways it can be written down as a linear combination with nonnegative integer coefficients of the positive roots. The Kostant partition function is a simple example of a more general class of functions, called \emph{vector partition functions}.

\begin{definition}
Let $M$ be a $d \times n$ matrix over the integers, such that $\mathrm{ker}M \cap \mathbb{R}^n_{\geq 0} = 0$. The \emph{vector partition function} (or simply \emph{partition function}) associated to $M$ is the function
\begin{displaymath}
\begin{array}{rccl}
\phi_M : & \mathbb{Z}^d & \longrightarrow & \mathbb{N}\\
& b & \mapsto & |\{x \in \mathbb{N}^n \ : \  Mx = b\}|
\end{array}
\end{displaymath}
\end{definition}

The condition $\mathrm{ker}M \cap \mathbb{R}^n_{\geq 0} = 0$ forces the set $\{x \in \mathbb{N}^n \ : \  Mx = b\}$ to have finite size, or equivalently, the set $\{x \in \mathbb{R}^n_{\geq 0} \ : \  Mx = b\} $ to be compact, in which case it is a polytope $P_b$, and the partition function is the number of integral points (lattice points) inside it.

Also, if we let $M_1, \ldots, M_n$ denote the columns of $M$ (as column-vectors), and $x=(x_1, \ldots, x_n) \in \mathbb{R}^n_{\geq 0}$, then $Mx = x_1M_1 + x_2M_2 + \cdots + x_nM_n$ and for this to be equal to $b$, $b$ has to lie in the cone $\mathrm{pos}(M)$ spanned by the vectors $M_i$. So $\phi_M$ vanishes outside of $\mathrm{pos}(M)$. 

It is well-known that partition functions are piecewise quasipolynomial, and that the domains of quasipolynomiality form a complex of convex polyhedral cones, called the \emph{chamber complex}. Sturmfels gives a very clear explanation in \cite{Sturmfels} of this phenomenon. The explicit description of the chamber complex is due to Alekseevskaya, Gel'fand and Zelevinski$\breve{\i}$ \cite{AGZ}. There is a special class of matrices for which partition functions take a much simpler form. Call an integer $d\times n$ matrix $M$ of full rank $d$ \emph{unimodular} if every nonsingular $d\times d$ submatrix has determinant $\pm 1$. For unimodular matrices, the chamber complex determines domains of polynomiality instead of quasipolynomiality \cite{Sturmfels}.

It is useful for what follows to describe how to obtain the chamber complex of a partition function. Let $M$ be a $d\times n$ integer matrix of full rank $d$ and $\phi_M$ its associated partition function. For any subset $\sigma \subseteq \{1,\ldots,n\}$, denote by $M_\sigma$ the submatrix of $M$ with column set $\sigma$, and let $\tau_\sigma = \mathrm{pos}(M_\sigma)$, the cone spanned by the columns of $M_\sigma$. Define the set $\mathcal{B}$ of \emph{bases} of $M$ to be
\begin{displaymath}
\mathcal{B} = \{\sigma \subseteq \{1,\ldots,n\}\ : \ |\sigma| = d \ \ \textrm{and } \ \mathrm{rank}(M_\sigma) = d\}\,.
\end{displaymath}
$\mathcal{B}$ indexes the invertible $d\times d$ submatrices of $M$. The \emph{chamber complex} of $\phi_M$ is the common refinement of all the cones $\tau_\sigma$, as $\sigma$ ranges over $\mathcal{B}$ (see \cite{AGZ}). A theorem of Sturmfels \cite{Sturmfels} describes exactly how partition functions are quasipolynomial over the chambers of that complex. 

If we let $M_{A_n}$ be the matrix whose columns are the positive roots $\Delta_+^{(A_n)}$ of $A_n$, written in the basis of simple roots, then we can write Kostant's partition function in the matrix form defined above as
\begin{displaymath}
K_{A_n}(v) = \phi_{M_{A_n}}(v)\,.
\end{displaymath}

The following lemma is a well-known fact about $M_{A_n}$ and can be deduced from general results on matrices with columns of $0$'s and $1$'s where the $1$'s come in a consecutive block (see \cite{Schrijver}).

\begin{lemma}
The matrix $M_{A_n}$ is unimodular for all $n$. 
\end{lemma}

$M_{A_n}$ unimodular means that the Kostant partition functions for $A_n$ is polynomial instead of quasipolynomial on the cells of the chamber complex. In general, for $M$ unimodular, the polynomial pieces have degree at most the number of columns of the matrix minus its rank (see \cite{Sturmfels}). In our case, $M_{A_n}$ has rank $n$ and as many columns as $A_n$ has positive roots, ${n+1\choose 2}$. Hence the Kostant partition function for $A_n$ is piecewise polynomial of degree at most ${n+1\choose 2}-n = {n\choose 2}$.

\begin{remark}
\label{rem:LocallyPolynomial}
In view of Steinberg's formula \eqref{eqn:Steinberg}, this means that the Littlewood-Richardson coefficients are given by a piecewise polynomial function of degree at most ${n\choose 2}$ in the three sets of variables $\lambda$, $\mu$ and $\nu$, if these partitions have at most $n+1$ parts. This will be made precise in Sections~\ref{sec:SteinbergArrangement} and~\ref{sec:PolynomialityComplex}
\end{remark}

\section{A vector partition function for the Littlewood-Richardson coefficients}
\label{sec:PartitionFunction}

There are many combinatorial ways to compute the Littlewood-Richardson coefficients, in particular the Littlewood-Richardson rule \cite{EC2}, honeycombs \cite{KnutsonTao} and Berenstein-Zelevinsky triangles \cite{BZ2}. The model that is most convenient for us is the hive model \cite{Buch,KnutsonTao}.

\begin{definition}
A \emph{$k$-hive} is an array of numbers $a_{ij}$ with $0\leq i,j\leq k$ and $i+j\leq k$. We will represent hives in matrix form. For example, a $4$-hive is
\begin{equation}
\renewcommand{\arraystretch}{1.8}
\begin{array}{lllll}
a_{00} & a_{01} & a_{02} & a_{03} & a_{04} \\
a_{10} & a_{11} & a_{12} & a_{13}\\
a_{20} & a_{21} & a_{22}\\
a_{30} & a_{31}\\
a_{40}
\end{array}
\renewcommand{\arraystretch}{1.0}
\end{equation}
We will call a hive \emph{integral} if all its entries are nonnegative integers
\end{definition}

Following the terminology of \cite{KTT}, we will call \emph{hive conditions} (HC) the conditions
\begin{equation}
\begin{array}{c@{\hspace{0.1\textwidth}}c@{\hspace{0.1\textwidth}}c}
\renewcommand{\arraystretch}{1.8}
\begin{array}{ccc} 
& j & j+1\\
i & \bullet & \fbox{$\bullet$}\\
i+1 & \fbox{$\bullet$} & \bullet
\end{array} & 
\renewcommand{\arraystretch}{1.8}
\begin{array}{ccc} 
& j & j+1\\
i & & \bullet\\
i+1 & \fbox{$\bullet$} & \fbox{$\bullet$}\\
i+2 & \bullet &
\end{array} & 
\renewcommand{\arraystretch}{1.8}
\begin{array}{cccc}
& \phantom{+}j\phantom{0} & j+1 & j+2 \\
i & & \fbox{$\bullet$} & \bullet\\
i+1 & \bullet & \fbox{$\bullet$} &
\end{array}
\end{array}
\end{equation}
where in each diagram, the sum of the boxed entries is at least as large as the sum of the other two entries. In terms of the $a_{ij}$, (HC) is
\begin{eqnarray} 
a_{i+1\,j}+a_{i\,j+1} & \geq & a_{ij}+a_{i+1\,j+1}\\
a_{i+1\,j}+a_{i+1\,j+1} & \geq & a_{i+2\,j}+a_{i\,j+1}\\
a_{i\,j+1}+a_{i+1\,j+1} & \geq & a_{i+1\,j}+a_{i\,j+2}
\end{eqnarray}
for $i+j\leq k-2$.

\begin{proposition}{\rm\textbf{(Knutson-Tao \cite{KnutsonTao}, Fulton \cite{Buch})}}
For $\lambda$, $\mu$ and $\nu$ partitions with at most $k$ parts and $|\lambda|+|\mu|=|\nu|$, the Littlewood-Richardson coefficient $c_{\lambda\mu}^{\nu}$ is the number of integral $k$-hives satisfying (HC) and the boundary conditions
\begin{equation}
\begin{array}{rcl@{\hspace{0.05\textwidth}}l}
a_{00} & = & 0\,, & \\[2mm]
a_{0j} & = & \lambda_1+\cdots+\lambda_j & 1\leq j\leq k\\[2mm]
a_{i0} & = & \nu_1+\cdots+\nu_i & 1\leq i\leq k\\[2mm]
a_{m,k-m} & = & |\lambda| + \mu_1+\cdots+\mu_m & 1\leq m\leq k\,.
\end{array}
\end{equation}
\end{proposition}

Once the boundary conditions are imposed, we are left with a system of inequalitites in the nonnegative integral variables $a_{ij}$ for $1\leq i,j\leq k-1$ and $i+j\leq k-1$. If we let these $a_{ij}$ take real values, the inequalitites define a rational polytope $Q_{\lambda\mu}^{\nu}$, and the Littlewood-Richardson coefficient corresponding to the boundary conditions is the number of integral (lattice) points inside $Q_{\lambda\mu}^{\nu}$. 

Given a $d$-dimensional rational polytope $Q$ in $\mathbb{R}^n$, we will denote by $mQ$ the polytope $Q$ blown up by a factor of $m$. The function $m\in\mathbb{N}\mapsto|mQ\cap\mathbb{Z}^n|$ is called the \emph{Ehrhart function} of $Q$, and is known \cite{Ehrhart,EC2} to be a quasipolynomial of degree $d$ in $m$. Furthermore, if $Q$ is integral, the Ehrhart function is a degree $d$ polynomial in $m$. This means that the function
\begin{equation}
N \longmapsto c_{N\lambda\,N\mu}^{N\nu}
\end{equation}
is the Ehrhart quasipolynomial of the polytope $Q_{\lambda\mu}^{\nu}$. It is known that $Q_{\lambda\mu}^{\nu}$ is not integral in general (see examples in \cite{KTT}).

This describes the behavior of the Littlewood-Richardson coefficients on a ray in $(\lambda,\mu,\nu)$-space, but we will get more general results by showing that we can find a vector partition function that gives these coefficients. We will then be able to work with conical chambers in $(\lambda,\mu,\nu)$-space instead of simple rays. This is accomplished in a way very similar to the one introduced for the weight multiplicities in \cite{BGR}, and this case is even simpler because the variables $a_{ij}$ are already constrained to be nonnegative.

We start by writing all the inequalities in the form
\begin{equation}
\sum_{{}\atop{\displaystyle1\leq i,j\leq k-1\atop {{}\atop \displaystyle i+j\leq k-1}}}b_{mij}a_{ij} \quad\leq\quad \sum_{1\leq t\leq k}c_{mt}\lambda_{t} + \sum_{1\leq t'\leq k}d_{mt'}\mu_{t'} + \sum_{1\leq t''\leq k}e_{mt''}\nu_{t''}\,,
\end{equation}
where $m$ indexes the inequalities. In a $k$-hive, there are ${k\choose 2}$ inequalities of the square type in the diagram above, and also ${k\choose 2}$ of them for each of the two parallelogram types. So we have $n(k) = 3{k\choose 2}$ inequalities overall and hence $1\leq m\leq n(k)$.

We next transform these inequalities into equalities by introducing a slack variable $s_m$ for each inequality:
\begin{equation}
\sum_{{}\atop{\displaystyle1\leq i,j\leq k-1\atop {{}\atop \displaystyle i+j\leq k-1}}}b_{mij}a_{ij} + s_m \quad=\quad  \sum_{1\leq t\leq k}c_{mt}\lambda_{t} + \sum_{1\leq t'\leq k}d_{mt'}\mu_{t'} + \sum_{1\leq t''\leq k}e_{mt''}\nu_{t''}\,.
\end{equation}

Solving the system of inequalities for nonnegative integral $a_{ij}$ is the same as solving the system of equalities for nonnegative integral $a_{ij}$ and $s_m$. Hence we are trying to solve the system
\begin{equation}
\underbrace{\begin{array}{@{}c@{}}\left(\begin{array}{ccc|ccc}
& & & & &\\
& b_{m,ij} & & & I_{n(k)} & \\
& & & & &
\end{array}\right)\\[-2mm]{}\end{array}}_{E_k}\cdot\left(\begin{array}{c}
a_{11}\\ \vdots\\ a_{1\,k-1}\\ \vdots\\ a_{k-1\,1}\\ s_1\\ \vdots\\ s_m
\end{array}\right) = \underbrace{\begin{array}{@{}c@{}}\left(\begin{array}{ccc|ccc|ccc}
& & & & & & & &\\
& c_{mt} & & & d_{mt'} & & & e_{mt''} &\\
& & & & & & & &
\end{array}\right)\\[-2mm]{}\end{array}}_{B_k}\cdot\left(\begin{array}{c}
\lambda_1\\ \vdots\\ \lambda_k\\ \mu_1\\ \vdots\\ \mu_k\\ \nu_1\\ \vdots\\ \nu_k\\
\end{array}\right)
\end{equation}
for integral nonnegative $a_{ij}$ and $s_m$. We have therefore proved the following.

\begin{theorem}
\label{thm:LRpartfunc}
The function $(\lambda,\mu,\nu)\mapsto c_{\lambda\mu}^{\nu}$ for $\lambda$, $\mu$, $\nu$ partitions with at most $k$ parts such that $|\lambda|+|\mu|=|\nu|$ and $\lambda,\mu\subseteq\nu$ is given by
\begin{equation}
c_{\lambda\mu}^{\nu} = \phi_{E_k}\left(B_k\left(\begin{array}{c}\lambda\\ \mu\\ \nu\end{array}\right)\right)\,.
\end{equation}
\end{theorem}

The chamber complex defined by $E_k$ is much too big for our purposes. For one thing, its cones have dimension $n(k)=3{k\choose 2}$, whereas $(\lambda,\mu,\nu)$-space is $3k$-dimensional. To simplify things, we can first restrict ourselves with the intersection of the complex of $E_k$ with the subspace
\begin{equation}
\mathcal{B}^{(k)} = \left\{\left(B_k\left(\begin{array}{c}\lambda\\ \mu\\ \nu\end{array}\right)\right)\ : \ \lambda,\mu,\nu\in\mathbb{R}^k\right\}
\end{equation}
of $\mathbb{R}^{n(k)}$ to get a complex $\mathcal{C}_k$. Then we can pull back the cones along the transformation $B_k$ to $(\lambda,\mu,\nu)$-space. Cones in $\mathcal{B}^{(k)}$ are given by inequalitites of the form
\begin{displaymath}
\left\langle v_i, B_k\left(\begin{array}{c}\lambda\\ \mu\\ \nu\end{array}\right)\right\rangle \geq 0
\end{displaymath}
for some directions $v_i\in\mathbb{R}^{n(k)}$. But
\begin{displaymath}
\left\langle v_i, B_k\left(\begin{array}{c}\lambda\\ \mu\\ \nu\end{array}\right)\right\rangle \geq 0 \quad\Leftrightarrow\quad \left\langle B_k^{\,T}v_i, \left(\begin{array}{c}\lambda\\ \mu\\ \nu\end{array}\right)\right\rangle \geq 0\,,
\end{displaymath}
where $B_k^{\,T}$ is the transpose of $B_k$. So we can pull back the cones to get a complex $B_k^{\,*}\mathcal{C}_k$ in $(\lambda,\mu,\nu)$-space. As a final simplification, we can note that $c_{\lambda\mu}^{\nu}=0$ unless $\lambda,\mu\subseteq\nu$ and $|\lambda|+|\mu|=|\nu|$ and that these conditions define a cone $C_k^{(1)}$ since the containment equations can be written $\lambda_i,\mu_i\leq\nu_i$ for $1\leq i\leq k$. The conditions $\lambda_1\geq\cdots\geq\lambda_k\geq 0$, $\mu_1\geq\cdots\geq\mu_k\geq 0$ and $\nu_1\geq\cdots\geq\nu_k\geq 0$ also define a cone $C_k^{(2)}$.

\begin{definition}
We will call the intersection of the cones $C_k^{(1)}$ and $C_k^{(2)}$ with the rectified complex $B_k^{\,*}\mathcal{C}_k$ the \emph{Littlewood-Richardson complex}, and denote it $\mathcal{LR}_k$. This complex lives on the subspace $|\lambda|+|\mu|=|\nu|$ of $\mathbb{R}^{3k}$.
\end{definition}

As a result of the general theory of vector partition functions, we get the following corollary.

\begin{corollary}
Under the conditions of the theorem above, the function $(\lambda,\mu,\nu)\mapsto c_{\lambda\mu}^{\nu}$ is quasipolynomial of degree at most $3{k\choose 2}+n(k)-\mathrm{rank}\,E_k = 3{k\choose 2}$ over the chambers of the complex $\mathcal{LR}_k$. 
\end{corollary}

We will show in Section~\ref{sec:PolynomialityComplex} that we actually get polynomials in the chambers.

It rapidly becomes computationally hard to work out the chamber complex and the associated polynomials; we present an example of how the computations are done on the simplest nontrivial example, $k=3$, in Section~\ref{sec:Keq3}.

\section{The Steinberg arrangement}
\label{sec:SteinbergArrangement}

In this section, we will construct a hyperplane arrangement whose regions are domains of polynomiality for the Littlewood-Richardson coefficients. We will deduce the form of this arrangement from a closer look at Steinberg's formula~(\ref{eqn:Steinberg}) and the chamber complex of the Kostant partition function defined in Section~\ref{subsec:ChamberComplexes}. 

The following lemma, proved in \cite{BGR} but reproduced here for the sake of completeness, describes the set of normals to the hyperplanes supporting the cells of the chamber complex for the Kostant partition function. 

\begin{lemma}
\label{lemma:KPFwalls}
The set of normals to the facets of the maximal cones of the chamber complex of the Kostant partition function of $A_n$ consists of all the conjugates of the fundamental weights.
\end{lemma}

\begin{proof}
The facets of the maximal cones of the chamber complex span the same hyperplanes as the facets of the base cones whose common refinement is the chamber complex. Base cones correspond to sets of $n$ linearly independent positive roots.  Fixing a particular base cone spanned by $\{\gamma_1,\ldots,\gamma_n\}$, consider the undirected graph $G$ on $\{1,\ldots,n+1\}$ where $(i,j)$ is an edge if $e_i-e_j=\gamma_m$ for some $m$. The fact that the $\gamma_j$'s are linearly independent implies that $G$ has no cycles. So $G$ is a forest, and since it has $n+1$ vertices and $n$ edges (one for each $\gamma_j$), it is actually a tree. Suppose now we remove $\gamma_j = e_s-e_t$ and want to find the normal of the hyperplane spanned by the other $\gamma_i$'s. The graph $G$ with the edge $(s,t)$ removed consists of two trees $T_1$ and $T_2$. List $\{1,\ldots,n+1\}$ in the form
\begin{displaymath}
\sigma : \ \underbrace{i_1, i_2, \ldots, i_{j-1}, s}_{\textrm{vertices of $T_1$}}, \underbrace{t, i_j, i_{j+1}, \ldots, i_{n+1-2}}_{\textrm{vertices of $T_2$}}
\end{displaymath}
where we will think of $\sigma$ as a permutation in one-line form. 

Now let $\alpha_i'=e_{\sigma(i)}-e_{\sigma(i+1)}$ and note that $\alpha_j'=e_s-e_t=\gamma_j$. The set $\{\alpha_1',\ldots,\alpha_n'\}$ is a root system basis because it is the image under the action of $\sigma^{-1}$ of the original simple roots $\alpha_i=e_i-e_{i+1}$. Observe that every edge in $T_1$ can be expressed as a sum of $\alpha_1',\ldots,\alpha_{j-1}'$, and every edge in $T_2$ as a sum of $\alpha_{j+1}',\ldots,\alpha_n'$, so that all $\gamma_i$'s in $\{\gamma_1,\ldots,\widehat{\gamma_j},\ldots,\gamma_n\}$ can be expressed as linear combinations of $\alpha_1',\ldots,\widehat{\alpha_j'},\ldots,\alpha_n'$. The normal for the corresponding hyperplane will therefore be the $j$th fundamental weight $\omega_j'$ for the basis $\{\alpha_1',\ldots,\alpha_n'\} = \sigma\cdot\{\alpha_1\ldots,\alpha_n\}$.

Conversely, given any fundamental weight $\omega_j'$ for the root system basis $\sigma\cdot\{\alpha_1\ldots,\alpha_n\}$ (or equivalently, $\sigma^{-1}\cdot\omega_j$, where $\omega_j$ is the $j$th fundamental weight for the standard simple roots), we want to show it can occur as the normal to a hyperplane. Let $H$ be a hyperplane separating the standard positive roots from the negative ones. For each $\alpha_i'=\sigma\cdot\alpha_i$, we can pick a sign $\varepsilon_i$ such that $\varepsilon_i\alpha_i'$ is on the positive side of $H$. Hence $\{\varepsilon_1\alpha_1',\ldots,\varepsilon_n\alpha_n'\}$ is a linearly independent subset of the set of standard positive roots, and thus it corresponds to one of the base cones of $M_{A_n}$. The corresponding graph is a path since we have a system of simple roots (up to sign reversal). Removing $\varepsilon_j\alpha_j'$ and applying the above procedure with the order given by the path gives that $\omega_j'$ occurs as the normal of the corresponding hyperplane.
\end{proof}

To compute the Littlewood-Richardson coefficients using Steinberg's formula~\eqref{eqn:Steinberg}, we look at the points $\sigma(\lambda+\delta)+\tau(\mu+\delta)-(\nu+2\delta)$, as $\sigma$ and $\tau$ range over the Weyl group $\mathfrak{S}_k$ (we assume here that $\lambda$, $\mu$ and $\nu$ have at most $k$ parts and index irreducible representations of $\gln{k}$). Some of these points will lie inside the chamber complex for the Kostant partition function and we compute the Littlewood-Richardson coefficients by finding which cells contain them and evaluating the corresponding polynomials at those points. We will call $(\lambda,\mu,\nu)$ \emph{generic} if none of the points $\sigma(\lambda+\delta)+\tau(\mu+\delta)-(\nu+2\delta)$ lies on a wall of the chamber complex of the Kostant partition function. If we change a generic $(\lambda,\mu,\nu)$ to $(\lambda',\mu',\nu')$ on the hyperplane $|\lambda|+|\mu|=|\nu|$ in such a way that none of the $\sigma(\lambda+\delta)+\tau(\mu+\delta)-(\nu+2\delta)$ crosses a wall, we will obtain $c_{\lambda'\mu'}^{\nu'}$ by evaluating the same polynomials. So there is a neighborhood of $(\lambda,\mu,\nu)$ on which the Littlewood-Richardson coefficients are given by the same polynomial in the variables $\lambda$, $\mu$ and $\nu$. 

Lemma~\ref{lemma:KPFwalls} describes the walls of the chamber complex for the Kostant partition function in terms of the normals to the hyperplanes (though the origin) supporting the facets of the maximal cells. Now a point $\sigma(\lambda+\delta)+\tau(\mu+\delta)-(\nu+2\delta)$ will be on one of those walls (hyperplane though the origin) when its scalar product with the hyperplane's normal, say $\theta(\omega_j)$, vanishes, that is when
\begin{equation}\label{eqn:KosArrWalls}
\langle \sigma(\lambda+\delta)+\tau(\mu+\delta)-(\nu+2\delta), \theta(\omega_j)\rangle = 0
\end{equation}

Consider the arrangement on the subspace $|\lambda|+|\mu|=|\nu|$ of $\mathbb{R}^{3k}$ consisting of all such hyperplanes, for $1\leq j\leq k$ and $\sigma,\,\tau,\,\theta\in\mathfrak{S}_k$. For $(\lambda,\mu,\nu)$ and $(\lambda',\mu',\nu')$ in the same region of this arrangement and any fixed $\sigma,\tau\in\mathfrak{S}_k$, the points $\sigma(\lambda+\delta)+\tau(\mu+\delta)-(\nu+2\delta)$ and $\sigma(\lambda'+\delta)+\tau(\mu'+\delta)-(\nu'+2\delta)$ lie on the same side of every wall of the chamber complex for the Kostant partition function. We will call this arrangement the \emph{Steinberg arrangement}, and denote it $\mathcal{SA}_k$.

\begin{definition}
Fix a labelling on the chambers of the complex for the Kostant partition function, and let $p_1,\,p_2,\,\ldots$ be the polynomials associated to the chambers. For generic $\lambda$, $\mu$ and $\nu$, let $v_{\sigma\tau}(\lambda,\mu,\nu)$ be the label of the region containing the point $\sigma(\lambda+\delta)+\tau(\mu+\delta)-(\nu+2\delta)$ (this label is unique for generic $\lambda$, $\mu$ and $\nu$). Define the \emph{type} of $\lambda$, $\mu$ and $\nu$ to be the matrix
\begin{displaymath}
\mathrm{Type}(\lambda,\mu,\nu) = \big(v_{\sigma\tau}(\lambda,\mu,\nu)\big)_{\sigma,\tau\in\mathfrak{S}_k}\,,
\end{displaymath}
for some fixed total order on $\mathfrak{S}_k$. Furthermore, define 
\begin{equation}
\label{eqn:TypePolynomial}
P(\lambda,\mu,\nu) = \sum_{\sigma\in\mathfrak{S}_k}\sum_{\tau\in\mathfrak{S}_k}(-1)^{\mathrm{inv}(\sigma\tau)}p_{v_{\sigma\tau}}(\sigma(\lambda+\delta)+\tau(\mu+\delta)-(\nu+2\delta))\,.
\end{equation}
\end{definition}

\begin{proposition}
\label{prop:SteinbergPolynomial}
$P(\lambda,\mu,\nu)$ is a polynomial function in $\lambda$, $\mu$ and $\nu$ on the interior of the regions of $\mathcal{SA}_k$ and gives the Littlewood-Richardson coefficients there.
\end{proposition}

\begin{proof}
The type of points along a path between $(\lambda',\mu',\nu')$ and $(\lambda'',\mu'',\nu'')$ in the interior of the same region of $\mathcal{SA}_k$ will remain the same by definition of the Steinberg arrangement (because no $\sigma(\lambda+\delta)+\tau(\mu+\delta)-(\nu+2\delta)$ crosses a wall along that path).
\end{proof}

The reason why Proposition~\ref{prop:SteinbergPolynomial} is restricted to the interior of the regions is that while polynomials for adjacent regions of the chamber complex for the Kostant partition function have to coincide on the intersection of their closures, there is a discontinuous jump in the value of the Kostant partition function (as a piecewise polynomial function) when going from a region on the boundary of the complex to region $0$ (outside the complex). 

To summarize, the hyperplanes of the Steinberg arrangement are defined by the equations 
\begin{equation}
\langle \sigma(\lambda+\delta)+\tau(\mu+\delta)-(\nu+2\delta), \theta(\omega_j)\rangle = 0 
\end{equation}
or
\begin{equation}
\label{eqn:SteinbergHyp}
\langle \sigma(\lambda)+\tau(\mu)-\nu, \theta(\omega_j)\rangle = \langle 2\delta-\sigma(\delta)-\tau(\delta), \theta(\omega_j)\rangle\,.
\end{equation}

Note that the right hand side of ~\eqref{eqn:SteinbergHyp} doesn't depend on $\lambda$, $\mu$ and $\nu$, and we will call it the $\emph{$\delta$-shift}$: 
\begin{equation}
s(\sigma, \tau, \theta, j) = \langle 2\delta-\sigma(\delta)-\tau(\delta), \theta(\omega_j)\rangle\,.
\end{equation}

\section{Polynomiality in the chamber complex}
\label{sec:PolynomialityComplex}

We have now expressed the Littlewood-Richardson coefficients in two ways: as a quasipolynomial function over the cones of the chamber complex $\mathcal{LR}_k$, and as a polynomial function over the interior of the regions of the hyperplane arrangement $\mathcal{SA}_k$. In this section, we relate the chamber complex to the hyperplane arrangement to show that the quasipolynomials are actually polynomials.

\begin{theorem}
\label{thm:PolynomialityComplex}
The quasipolynomials giving the Littlewood-Richardson coefficients in the cones of the chamber complex $\mathcal{LR}_{k}$ are polynomials of degree at most ${k-1 \choose 2}$ in the three sets of variables $\lambda=(\lambda_1,\ldots,\lambda_k)$, $\mu=(\mu_1,\ldots,\mu_k)$ and $\nu=(\nu_1,\ldots,\nu_k)$. 
\end{theorem}

\begin{proof}
We will show that for each cone $C$ of $\mathcal{LR}_{k}$ we can find a region $R$ of the Steinberg arrangement $\mathcal{SA}_k$ such that $C\cap R$ contains an arbitrarily large ball. Then $P(\lambda,\mu,\nu)$ and the quasipolynomial in $C$ agree on the lattice points $(\lambda,\mu,\nu)$ in that ball, and must therefore be equal. The degree bounds follow from the degree bounds on the polynomials giving the Kostant partition function (see Remark~\ref{rem:LocallyPolynomial}). Note that since $c_{\lambda\mu}^{\nu}$ is invariant under adding ``0'' parts to the partitions, we get the best degree bound by working in $\sln{k}$ for $k$ as small as possible, that is $k=\max\{l(\lambda),l(\mu),l(\nu)\}$. 

We can deform $\mathcal{SA}_k$ continuously to make the $\delta$-shifts zero, by considering the arrangement $\mathcal{SA}_k^{(t)}$ with hyperplanes
\begin{equation}
\langle \sigma(\lambda)+\tau(\mu)-\nu, \theta(\omega_j)\rangle = t\,\langle 2\delta-\sigma(\delta)-\tau(\delta), \theta(\omega_j)\rangle
\end{equation}
and letting $t$ going from 1 to 0, for example. The final deformed arrangement $\mathcal{SA}_k^{(0)}$ is a central arrangement (all the hyperplanes go through the origin) whose regions are therefore cones. $C$ will intersect nontrivially one of the cones $\tilde{R}$ of this arrangement (i.e. the dimension of the cone $C\cap\tilde{R}$ is the same as that of $C$ and $\tilde{R}$). Let $R$ be any region of $\mathcal{SA}_k$ whose deformed final version is $\tilde{R}$. Consider a ball of radius $r$ inside $\tilde{R}\cap C$, and suppose it is centered at the point $x$. Let $s$ is the maximal amount by which the hyperplanes of the Steinberg arrangement are shifted, i.e.
\begin{equation}
s = \max_{{}\atop{\displaystyle\sigma,\tau,\theta\in\mathfrak{S}_k\atop \displaystyle 1\leq j\leq k}}|\langle 2\delta-\sigma(\delta)-\tau(\delta), \theta(\omega_j)\rangle|\,.
\end{equation}
Then $R$ contains the ball of radius $r-s$ centered at $x$, and so does $C\cap R$. Since $C$ is a cone, we can make $r$ arbitrary large and the result follows since $s$ is bounded for fixed $k$.
\end{proof}

From this, we can deduce a ``stretching'' property for Littlewood-Richardson coefficients.

\begin{corollary}
\label{cor:Stretching}
The Littlewood-Richardson coefficients $c_{N\lambda\,N\mu}^{N\nu}$ are given by a polynomial in $N$ with rational coefficients. This polynomial has degree at most $3{k-1\choose 2}$ in $N$.
\end{corollary}

\begin{remark}
King, Tollu and Toumazet conjectured in \cite{KTT} that the $c_{N\lambda\,N\mu}^{N\nu}$ are polynomial in $N$ with nonnegative rational coefficients (Conjecture~\ref{conj:KTT} above). Corollary~\ref{cor:Stretching} establishes this conjecture, except for the nonnegativity of the coefficients. Derksen and Weyman \cite{DerksenWeyman} have a proof of this part of the conjecture using semi-invariants of quivers, and Knutson \cite{DerksenWeyman,Knutson} a proof using symplectic geometry techniques. 
\end{remark}

In fact, we can prove something stronger: we can perturb $(\lambda,\mu,\nu)$ a bit and get a more global stretching property.

\begin{corollary}
\label{cor:StableStretching}
Let $\Upsilon$ be the set 
\begin{equation}
\Upsilon = \{(\lambda,\mu,\nu)\ : \ \max\{l(\lambda),l(\mu),l(\nu)\}\leq k, |\lambda|+|\mu|=|\nu|, \lambda,\mu\subseteq\nu\}. 
\end{equation}
For any generic $(\lambda,\mu,\nu)\in\Upsilon$ we can find a neighborhood $U$\! of that point over which the function 
\begin{equation}
(\lambda,\mu,\nu,t)\in (U\cap\Upsilon)\times\mathbb{N} \longmapsto c_{t\lambda\,t\mu}^{t\nu}
\end{equation}
is polynomial of degree at most $3{k-1\choose 2}$ in $t$ and ${k-1\choose 2}$ in the $\lambda$, $\mu$ and $\nu$ coordinates.
\end{corollary}

\begin{proof}
Let $(\lambda,\mu,\nu) \in \Upsilon$. For $U$ sufficiently small, the points $\{(t\lambda,t\mu,t\nu)\ : \ t\in\mathbb{N}\}$ lie in the same cone of the chamber complex $\mathcal{LR}_{k}$. Hence the corresponding Littlewood-Richardson coefficients are obtained by evaluating the same polynomial at those points.  
\end{proof}

\section{An example for partitions with at most 3 parts}
\label{sec:Keq3}

We want to find a vector partition function counting the number of integral $3$-hives of the form
\begin{equation}
\renewcommand{\arraystretch}{2.5}
\begin{array}{llll}
0 & \lambda_1 & \lambda_1 + \lambda_2 & |\lambda|  \\
\nu_1 & a_{11} & |\lambda|+\mu_1 \\
\nu_1+\nu_2 & |\nu|-\mu_3 \\
|\nu|
\end{array}
\end{equation}

The hives conditions are given by
\begin{small}
\begin{equation}
\begin{array}{r@{\,\,}c@{\,\,}lr@{\,\,}c@{\,\,}lr@{\,\,}c@{\,\,}l}
a_{11} & \leq & \nu_1+\lambda_1 & -a_{11} & \leq & -\lambda_2-\nu_1 & -a_{11} & \leq & -\lambda_1-\nu_2 \\
-a_{11} & \leq & -\lambda_1-\lambda_3-\mu_1 & a_{11} & \leq & \lambda_1+\lambda_2+\mu_1 & -a_{11} & \leq & -\lambda_1-\lambda_2-\mu_2 \\
-a_{11} & \leq & -\lambda_1-\lambda_2-\lambda_3-\mu_1-\mu_2+\nu_2 & -a_{11} & \leq & \mu_2-\nu_1-\nu_2 & a_{11} & \leq & \lambda_1+\lambda_2+\lambda_3+\mu_1+\mu_2-\nu_3 
\end{array}
\end{equation}
\end{small}

This corresponds to the matrix system
\begin{footnotesize}
\begin{equation}
\underbrace{\begin{array}{@{}c@{}}\left(\begin{array}{rrrrrrrrrr}
1  & 1 & 0 & 0 & 0 & 0 & 0 & 0 & 0 & 0 \\
-1 & 0 & 1 & 0 & 0 & 0 & 0 & 0 & 0 & 0 \\
-1 & 0 & 0 & 1 & 0 & 0 & 0 & 0 & 0 & 0 \\
-1 & 0 & 0 & 0 & 1 & 0 & 0 & 0 & 0 & 0 \\
1  & 0 & 0 & 0 & 0 & 1 & 0 & 0 & 0 & 0 \\
-1 & 0 & 0 & 0 & 0 & 0 & 1 & 0 & 0 & 0 \\
-1 & 0 & 0 & 0 & 0 & 0 & 0 & 1 & 0 & 0 \\
-1 & 0 & 0 & 0 & 0 & 0 & 0 & 0 & 1 & 0 \\
1  & 0 & 0 & 0 & 0 & 0 & 0 & 0 & 0 & 1 \\
\end{array}\right)\\[-2mm]{}\end{array}}_{E_3}\cdot\left(\begin{array}{c}
a_{11} \\ s_1 \\ s_2 \\ \vdots \\ s_9
\end{array}\right) = \underbrace{\begin{array}{@{}c@{}}\left(\begin{array}{rrrrrrrrr}
 1 &  0 &  0     &  0 &  0 &  \phantom{+}0     &  1 &  0 &  0 \\
-1 &  0 & -1     & -1 &  0 &  0     &  0 &  0 &  0 \\
-1 & -1 & -1     & -1 & -1 &  0     &  0 &  1 &  0 \\
 0 & -1 &  0     &  0 &  0 &  0     & -1 &  0 &  0 \\
 1 &  1 &  0     &  1 &  0 &  0     &  0 &  0 &  0 \\
 0 &  0 &  0     &  0 &  1 &  0     & -1 & -1 &  0 \\
-1 &  0 &  0     &  0 &  0 &  0     &  0 & -1 &  0 \\
-1 & -1 &  0     &  0 & -1 &  0     &  0 &  0 &  0 \\
 1 &  1 &  1     &  1 &  1 &  0     &  0 &  0 & -1 
\end{array}\right)\\[-2mm]{}\end{array}}_{B_3}\cdot\left(\begin{array}{c}
\lambda_1 \\ \lambda_2 \\ \lambda_3 \\ \mu_1 \\ \mu_2 \\ \mu_3 \\ \nu_1 \\ \nu_2 \\ \nu_3
\end{array}\right)
\end{equation}
\end{footnotesize}

Note that $\mu_2$ doesn't not appear in this system. This is because it is determined by $|\lambda|+|\mu|=|\nu|$; we could have chosen another variable to disappear.

To get the chamber complex for the vector partition function associated to $E_3$, we have to find the sets of columns determining maximal nonsingular square matrices in $E_3$. These determine the bases cones whose common refinement gives the chamber complex. In our case, all subsets of $9$ columns determine a nonsingular matrix, so we get $10$ base cones. We can find their common refinement using a symbolic calculator like Maple or Mathematica; here we used Maple (version 8) and the package $\texttt{convex}$ by Matthias Franz \cite{convex}. We find the chamber complex $\mathcal{LR}_3$ by rectifying the cones to $(\lambda,\mu,\nu)$-space using $B_3^{\,T}$ and intersecting them with the cones $C_3^{(1)}$ and $C_3^{(2)}$. The list of rays of the cones of $\mathcal{LR}_3$
\begin{displaymath}
\begin{array}{c}
a_1 = (\,1\ 1\ 1\ |\ 0\ 0\ 0\ |\ 1\ 1\ 1\,) \qquad\qquad a_2 = (\,0\ 0\ 0\ |\ 1\ 1\ 1\ |\ 1\ 1\ 1\,) \\[2mm]
b =   (\,2\ 1\ 0\ |\ 2\ 1\ 0\ |\ 3\ 2\ 1\,) \\[2mm]
c = (\,1\ 1\ 0\ |\ 1\ 1\ 0\ |\ 2\ 1\ 1\,) \\[2mm]
d_1 = (\,1\ 1\ 0\ |\ 1\ 0\ 0\ |\ 1\ 1\ 1\,) \qquad\qquad d_2 = (\,1\ 0\ 0\ |\ 1\ 1\ 0\ |\ 1\ 1\ 1\,) \\[2mm]
e_1 = (\,1\ 1\ 0\ |\ 0\ 0\ 0\ |\ 1\ 1\ 0\,) \qquad\qquad e_2 = (\,0\ 0\ 0\ |\ 1\ 1\ 0\ |\ 1\ 1\ 0\,) \\[2mm]
f = (\,1\ 0\ 0\ |\ 1\ 0\ 0\ |\ 1\ 1\ 0\,) \\[2mm]
g_1 = (\,1\ 0\ 0\ |\ 0\ 0\ 0\ |\ 1\ 0\ 0\,) \qquad\qquad g_2 = (\,0\ 0\ 0\ |\ 1\ 0\ 0\ |\ 1\ 0\ 0\,) 
\end{array}
\end{displaymath}
where the bars separate the entries corresponding to the sets of variables $\lambda$, $\mu$ and $\nu$.

The following table gives the maximal ($8$-dimensional) cones of $\mathcal{LR}_3$, as well as the polynomial associated to each (computed by polynomial interpolation).

\begin{displaymath}
\renewcommand{\arraystretch}{1.5}
\begin{array}{|@{\hspace{5mm}}c@{\hspace{5mm}}|@{\hspace{5mm}}l@{\hspace{5mm}}|@{\hspace{5mm}}l@{\hspace{5mm}}|}
\hline
\textbf{Cone} & \textbf{Positive hull description} & \textbf{Polynomial} \\
\hline
\hline
\kappa_{1} & \mathrm{pos}(a_1, a_2, b, c, d_1, d_2, e_1, e_2) & 1 - \lambda_2 - \mu_2 + \nu_1 \\[1mm]
\hline
\kappa_{2} & \mathrm{pos}(a_1, a_2, b, c, d_1, d_2, g_1, g_2) & 1 + \nu_2 - \nu_3 \\[1mm]
\hline
\kappa_{3} & \mathrm{pos}(a_1, a_2, b, c, e_1, e_2, g_1, g_2) & 1 + \lambda_1 + \mu_1 - \nu_1 \\[1mm]
\hline
\kappa_{4} & \mathrm{pos}(a_1, a_2, b, d_1, d_2, e_1, e_2, f) & 1 + \nu_1 - \nu_2 \\[1mm]
\hline
\kappa_{5} & \mathrm{pos}(a_1, a_2, b, d_1, d_2, f, g_1, g_2) & 1 + \lambda_2 + \mu_2 - \nu_3 \\[1mm]
\hline
\kappa_{6} & \mathrm{pos}(a_1, a_2, b, e_1, e_2, f, g_1, g_2) & 1 - \lambda_3 - \mu_3 + \nu_3 \\[1mm]
\hline
\kappa_{7} & \mathrm{pos}(a_1, a_2, b, c, d_1, d_2, e_1, g_1) & 1 + \lambda_3 + \mu_1 - \nu_3 \\
\kappa_{8} & \mathrm{pos}(a_1, a_2, b, c, d_1, d_2, e_2, g_2) & 1 + \lambda_1 + \mu_3 - \nu_3 \\[1mm]
\hline
\kappa_{9\ } & \mathrm{pos}(a_1, a_2, b, c, d_1, e_1, e_2, g_2) & 1 + \lambda_1 - \lambda_2 \\
\kappa_{10} & \mathrm{pos}(a_1, a_2, b, c, d_2, e_1, e_2, g_1) & 1 + \mu_1 - \mu_2 \\[1mm]
\hline
\kappa_{11} & \mathrm{pos}(a_1, a_2, b, c, d_1, e_1, g_1, g_2) & 1 - \lambda_2 - \mu_3 + \nu_2 \\
\kappa_{12} & \mathrm{pos}(a_1, a_2, b, c, d_2, e_2, g_1, g_2) & 1 - \lambda_3 - \mu_2 + \nu_2 \\[1mm]
\hline
\kappa_{13} & \mathrm{pos}(a_1, a_2, b, d_1, d_2, e_1, f, g_1) & 1 - \lambda_1 - \mu_3 + \nu_3 \\
\kappa_{14} & \mathrm{pos}(a_1, a_2, b, d_1, d_2, e_2, f, g_2) & 1 - \lambda_3 - \mu_1 + \nu_3 \\[1mm]
\hline
\kappa_{15} & \mathrm{pos}(a_1, a_2, b, d_1, e_1, f, g_1, g_2) & 1 + \mu_2 - \mu_3 \\
\kappa_{16} & \mathrm{pos}(a_1, a_2, b, d_2, e_2, f, g_1, g_2) & 1 + \lambda_2 - \lambda_3 \\[1mm]
\hline
\kappa_{17} & \mathrm{pos}(a_1, a_2, b, d_1, e_1, e_2, f, g_2) & 1 + \lambda_1 + \mu_2 - \nu_2 \\
\kappa_{18} & \mathrm{pos}(a_1, a_2, b, d_2, e_1, e_2, f, g_1) & 1 + \lambda_2 + \mu_1 - \nu_2 \\[1mm] 
\hline
\end{array}
\end{displaymath}

\begin{remark}
The symmetry $c_{\lambda\mu}^{\nu}=c_{\mu\lambda}^{\nu}$ implies that we can interchange the $\lambda$ and $\mu$ coordinates. This corresponds to a symmetry of the chamber complex $\mathcal{LR}_3$ under this transformation. This is why some of the rays and cones have been grouped in pairs. 
\end{remark}

\begin{remark}
We observe from the form of the polynomials in the table above that the equation 
\begin{equation}
c_{N\lambda\,N\mu}^{N\nu} = 1 + N(c_{\lambda\mu}^{\nu}-1)
\end{equation}
holds for $l(\lambda),l(\mu),l(\nu)\leq 3$. This was previously observed in \cite{KTT}.
\end{remark}

\bigskip

\begin{center}
{\large \textbf{Acknowledgments}}
\end{center}
\medskip

I would like to thank Sara Billey, Benoit Charbonneau, Victor Guillemin for useful discussions and comments, Ron King especially for telling me about his interesting conjectures, and also Allen Knutson and Anders Buch for bringing the results of \cite{DerksenWeyman} and \cite{Knutson} to my attention.

\medskip

\bibliographystyle{plain}

\bibliography{biblio}

\end{document}